\documentclass[12pt]{article}

\usepackage[numbered,framed]{matlab-prettifier}

\usepackage{enumitem}
\usepackage{listings}
\usepackage{adjustbox}
\usepackage[margin=1in,footskip=0.25in]{geometry}
\usepackage{relsize}
\usepackage{amsthm}%
\usepackage{mathtools}
\usepackage{amsfonts}%
\usepackage{pdfpages}
\usepackage{amsmath}
\usepackage{color}
\usepackage{flexisym}
\definecolor{mygreen}{RGB}{28,172,0} 
\definecolor{mylilas}{RGB}{170,55,241}

\usepackage{cite}
\newtheorem{theorem}{Theorem}
\newtheorem{remark}[theorem]{Remark}

\begin{document}

\title{Symmetry-breaking and bifurcation diagrams of fractional-order maps}

\author{Marius-F. Danca{\footnote{Corresponding author}}\\
Romanian Institute of Science and Technology, \\
400487, Cluj-Napoca, Romania\\
and\\
STAR-UBB Institute\\
Babes-Bolyai University,\\
400084, Cluj-Napoca, Romania\\
Email: danca@rist.ro\\
}

\maketitle

\begin{abstract}
In this paper two important aspects related to Caputo's fractional-order discrete variant of a class of maps defined on the complex plane, are analytically and numerically revealed: attractors symmetry-broken induced by the fractional-order and the sensible problem of determining the right bifurcation diagram of discrete systems of fractional-order. It is proved that maps of integer order with dihedral symmetry or cycle symmetry loose their symmetry once they are transformed in fractional-order maps. Also, it is conjectured that, contrarily to integer-order maps, determining the bifurcation diagrams of fractional-order maps is far from being a clarified problem. Two examples are considered: dihedral logistic map and cyclic logistic map.
\end{abstract}

\textbf{keyword }Symmetry breaking; Caputo fractional differences; FO logistic map; Bifurcation diagram;

\vspace{3mm}

\section{Introduction}
The whole of theoretical physics, and our general picture of the world, are based on symmetries \cite{first}. Symmetry breaking is a fundamental concept in particle physics. Thus, the most popular mechanism is the Higgs mechanism which has a key role in the electroweak theory, which unifies interactions via the weak force and the electromagnetic force \cite{brit}. As known, the logistic map of Integer Order (IO) presents a chaotic dynamic. The question is if, beside chaotic evolution this famous map, introduced by May \cite{uuu}, could also have some symmetry? For this purpose, the map should be modified. If one intends to obtain a symmetry on a line, consider the reflection (mirror or flip) symmetry, S, with respect the origin, which can be described as $S(x)=-x$. Considering an one-dimensional map $f:\mathbb{R}\rightarrow \mathbb{R}$, the symmetry about the origin means that considering two symmetric points $x$ and $S(x)$, the first iteration of these two points, i.e. $f(x)$ and $f(S(x))$, should remain symmetric with respect the origin, which means $f(S(x))=S(f(x))$. Because $S(x)=-x$, the above identity is verified if $f(-x)=-f(x)$, i.e. $f$ is an odd map.
Coming back to the logistic map, $f(x)=rx(1-x)$, which is not an odd map, in order to obtain symmetry, an odd modified variant $f(x)=rx(1-x^2)$ has been introduced by Field and Golubitsky in 1992 \cite{aa1}. There are many works on odd variants of the logistic map, one of the first papers being \cite{testa}.

Symmetries are easier to study in the complex plane than in the cartesian plane since, for some $z\in \mathbb{C}$, the reflection symmetry becomes $S(z)=\bar{z}$ and the rotation has the expression $R(z)=wz$, where $w$ is a complex number with unity modulus.

One of 2D symmetries is the dihedral symmetry presented by the \textbf{dihedral group} $D_m$, with $m=1,2,3,...$, one of the two series of discrete point groups in two dimensions.
In geometry, a two-dimensional point group or rosette group is a group of geometric symmetries (isometries) that keep at least one point fixed in a plane. In plane, $D_m$, a group with $2m$ elements, consists of $m$ rotations $R$ of multiples of $360^\circ/m$ about the origin, and $m$ symmetries $S$ across $m$ lines through the origin, making angles of multiples of $180^\circ/m$ with each other. This is the symmetry group of a regular polygon with $m$ sides ($m$-gon).
A figure with rotational symmetry appears the same after rotating by some amount around the center point.
Using degrees to describe the rotation amount is inconvenient because the precise angle is not obvious from looking at the figure. Therefore it is useful to use the order of rotation to describe rotational symmetry, i.e. a $\frac{2\pi}{m}$ radians angle of rotation about the center of the $m$-gon which leaves the figure unchanged. For example, the Mercedes Benz logo has dihedral $D_3$ symmetries.

Another spectacular group from the geometric point of view in the plane, considered in this paper, is the \textbf{cyclic group}, $C_m$ (a subgroup of $D_m$) where the only symmetries are rotations with $\frac{2\pi}{m}$.

Field and Golubitsky illustrated how to create chaotic attractors with both kind of symmetries in \cite{aa1}. Other recent work on creating chaotic attractors with symmetry \cite{aa2} has been
based on demanding that the coefficients of general polynomial or trigonometric functions
satisfy certain criterion in order to respect the symmetry group of interest. This has included
chaotic attractors with cubic \cite{treix}, tetrahedral \cite{doispex}, hypercube \cite{unspe}, frieze and wallpaper symmetries \cite{patrux}. A comprehensive literature on chaotic attractors with symmetry is presented in \cite{aa2}

An example of a map on the complex plane, $f:\mathbb{C}\rightarrow \mathbb{C}$ with dihedral or cyclic symmetry considered in this paper is \cite{aa1,chos} (see also \cite{dudu})

\begin{equation}\label{d1}
f(z)=(a+bz\bar{z}+ci)z+d \bar{z}^{\; m-1},~~m=1,2,3,...
\end{equation}

\noindent where $a,b,c,d,e$ are real parameters, $i^2=-1$ and, depending on the parameters values, $m$ indicates the dihedral or cyclic symmetry type and the number of the symmetry axes.

Discrete fractional difference equations represent nowadays
a new area for scientists \cite{xxx6},\cite{xxx2},\cite{cinci},\cite{xxx1},\cite{xxx8},\cite{opt},\cite{Go},\cite{04},\cite{doi},\cite{doispe}. Problems of discrete systems of FO, such as
hidden attractors and chaos control are analyzed in, e.g., \cite{hid},\cite{v1},\cite{v2},\cite{v3}. The extremely important problem of non-periodicity in Fractional Order (FO), continuous-time or discrete systems, are revealed in \cite{neper1} and \cite{nepe}, respectively.

In this paper the symmetry of the IO and Fractional Order (FO) variants of the logistic map \eqref{d1} are considered.

The paper is organized as follows: In Section 2 the symmetry of the IO variant of the map \eqref{d1} is considered; Section 3 deals with the lost of the symmetry of the map \eqref{d1} once the map is transformed into a FO map; Section 4 presents comparatively the bifurcation diagrams of the IO and FO variants of the map \eqref{d1} and the last section, Conclusion, summarizes the results.

\section{Symmetry in logistic maps \eqref{d1} of IO}
\subsection{Dihedral logistic map of IO}\label{sub21}

Let the map \eqref{d1} with $c=0$
\begin{equation}\label{d11}
f(z)=(a+bz\bar{z})z+d \bar{z}^{\;m-1}.
\end{equation}
The map is the \textbf{logistic map with dihedral symmetries}.

To prove that the logistic map \eqref{d11} is dihedral, consider the $m$-gon centered at the origin. Then, the elements of the $D_m$ dihedral group act as $2m$ linear transformations of the plane, based on the following matrix representations
\begin{equation*}
R_1=
\begin{pmatrix}
\cos\dfrac{2\pi}{m}& -\sin\dfrac{2\pi}{m}\\[12pt]
\sin\dfrac{2\pi}{m} &\cos\dfrac{2\pi}{m}
\end{pmatrix},
\end{equation*}
and
\begin{equation*}
S_0=
\begin{pmatrix}
1& 0\\
0&-1
\end{pmatrix}.
\end{equation*}
\noindent The other elements of the group $D_m$ are: $R_1^k$, and $S_k=R_kS_0$, for $k=0,1,2,...,m-1$, with the mention that $R_0=R_1^0=I_2$.


\noindent Beside the matrix representation of the $D_m$ elements, there exists the possibility to use the polar representation. Thus, if a complex number has the rectangular form, $z=x+iy$, the polar form is $z=r\mathrm{e}^{i\theta}$, where $r=\sqrt{x^2+y^2}$ is the modulus of $z$, while $\theta=\tan^{-1} (\frac{y}{x})$ is the argument (computationally determined with function $atan2(y,x)$). Therefore, by Euler's formula, the rotation and the reflection have the form

\begin{equation}\label{simus}
R_1(z)=r\mathrm{e}^{i(\theta+\frac{2\pi}{m})},
\end{equation}
and

\begin{equation*}\label{rota}
S_0(z)=\bar{z}=r\mathrm{e}^{-i\theta}.
\end{equation*}

\noindent Note that by rotation and reflection of a complex number, in this paper one understands that the transformations act on the image of the complex number.

\begin{theorem}\label{prop1}
The map $f$ defined by \eqref{d11}, presents dihedral symmetry.
\begin{proof}
\mbox{}
A function $f:\mathbb{C}\rightarrow \mathbb{C}$ is equivariant (concomitant) with respect to a symmetry $\sigma$ in $\mathbb{C}$ if \cite{aa2}
\begin{equation*}\label{sigma}
f(\sigma(z))=\sigma(f(z)), \text{~for all~} z\in \mathbb{C}.
\end{equation*}
Equivariant functions are a good choice for finding attractors with some desired symmetry. In the case of dihedral symmetries $\sigma$ is the rotation $R_1$ and the reflection $S_0$.

\noindent Let $z$ expressed in polar coordinates: $z=re^\theta$.
\begin{enumerate}
[leftmargin=0pt,align=left,itemindent=*,noitemsep]

\item
\emph{Reflection symmetry
}

Considering the rotation \eqref{simus} and the fact that $\bar{\bar{z}}=z$, $\bar{z}^k=\overline{z^k}$, and for $k\in \mathbb{R}$, $S_0(kz)=k\bar{z}$, one obtains
\begin{equation*}
f(S_0(z))=aS_0(z)+bS_0(z)S_0(\bar{z})S_0(z)+dS_0(\bar{z}^{\;m-1})=
a\bar{z}+b\bar{z}^2z+dz^{m-1}.
\end{equation*}
On the other side
\[
S_0(f(z))=\overline{az+bz^2 \bar{z}+e\bar{z}^{\;m-1}}=
a\bar{z}+b\bar{z}^2z+dz^{m-1},
\]
and, therefore $f(S_0(z))=S_0(f(z))$.
\item \emph{Rotation symmetry}
\[
f(R_1(z))=aR_1(z)+bR_1(z)\overline{R_1(z)}R_1(z)+d\overline{R_1(z)}^{\;m-1},
\]
and because $\overline{R_1(z)}=r\mathrm{e}^{-i[\theta+\frac{2\pi}{m}]}$ and $\overline{R_1(z)}^{\;m-1}=r^{m-1}\mathrm{e}^{-i(m-1)(\theta+\frac{2\pi}{m})}$, one obtains
\begin{equation*}
\begin{array}
[c]{ll}
f(R(z))=a  r \mathrm{e}^{i(\theta+\frac{2\pi}{m})}+b r\mathrm{e}^{i(\theta+\frac{2\pi}{m})}r\mathrm{e}^{-i(\theta+\frac{2\pi}{m})}r\mathrm{e}^{i(\theta+\frac{2\pi}{m})}+dr^{m-1}\mathrm{e}^{-i(m-1)(\theta+\frac{2\pi}{m})}\\[6pt]
=ar\mathrm{e}^{i(\theta+\frac{2\pi}{m})}+br^3\mathrm{e}^{i(\theta+\frac{2\pi}{m})}+dr^{m-1}\mathrm{e}^{-i(m-1)(\theta+\frac{2\pi}{m})}.
\end{array}
\end{equation*}
On the other side
\begin{equation*}
\begin{array}
[c]{cl}
R_1(f(z))&=R_1(az)+R_1(bz^2\bar{z})+R_1(d\bar{z}^{\;m-1})\\[6pt]
&=ar\mathrm{e}^{i(\theta+\frac{2\pi}{m})}+R_1(br^2\mathrm{e}^{2\theta}r\mathrm{e}^{-i\theta})+dR_1(\bar{z}^{\;m-1})\\[6pt]
&=ar\mathrm{e}^{i(\theta+\frac{2\pi}{m})}+R_1(br^3\mathrm{e}^{i\theta})+dr^{m-1}\mathrm{e}^{i\big[-(m-1)\theta+\frac{2\pi}{m}\big]}\\[6pt]
&=ar\mathrm{e}^{i(\theta+\frac{2\pi}{m})}+br^3\mathrm{e}^{i(\theta+\frac{2\pi}{m})}+dr^{m-1}\mathrm{e}^{i\big[-(m-1)\theta+\frac{2\pi}{m}\big]}.
\end{array}
\end{equation*}
and because $\frac{2\pi}{m}=-(m-1)\frac{2\pi}{m}+2\pi$, from the $2\pi$ periodicity it follows that $\mathrm{e}^{i[-(m-1)\theta+\frac{2\pi}{m}]}=\mathrm{e}^{i\big[-(m-1)(\theta+\frac{2\pi}{m})\big]}$, i.e. $f(R_1(z))=R_1(f(z))$,
and the theorem is proved.
 \qedhere
\end{enumerate}
\end{proof}
\end{theorem}

\noindent Similar results are obtained for the other elements $R_k$ and $S_k$.

\noindent To illustrate graphically the dihedral symmetry generated by the map \eqref{d1}, consider the following recursive iteration

\begin{equation}\label{d2}
z(n)=f(z(n-1)), ~~z(0)=z_0,~~n\in \mathbb{N}_1,
\end{equation}
with $\mathbb{N}_1=\{1,2,...\}$ and $z,z_0\in \mathbb{C}$, which is an Initial Value Problem (IVP) in the complex plane.

\noindent In this paper, in order to facilitate the numerical implementations, consider the notations $f(x,y):=(f_1(x,y),f_2(x,y))^T$, with $f_1(x,y)=\Re\{f(z)\}$ and $f_2(x,y)=\Im\{f(z)\}$, with which the IVP \eqref{d2} can be described in the rectangular plane as follows

\begin{equation}\label{unu}
\begin{array}
[c]{ll}%
\begin{pmatrix}
x_n\\
y_n
\end{pmatrix}
=
\begin{pmatrix}
f_1(x(n-1),y(n-1))\\[4pt]
f_2(x(n-1),y(n-1))
\end{pmatrix},\\[13pt]
~~ x_0,y_0\in \mathbb{R},~ n\in \mathbb{N}_1.
\end{array}
\end{equation}



Let $m=3$, for which one obtains the group $D_3$, i.e. the symmetries of the equilateral triangle. For $a=-1.804$, $b=1$, and $d=0.5$, the computational form of the IVP \eqref{unu} is rewritten in the following form

\begin{equation*}\label{xx}
\begin{array}
[c]{ll}%
\begin{pmatrix}
x_n\\
y_n
\end{pmatrix}
=
\begin{pmatrix}
-1.804~x_{n-1}+x_{n-1}^3+x_{n-1}y_{n-1}^2+0.5~x_{n-1}^2-0.5~y_{n-1}^2\\[4pt]
-1.804~y_{n-1}+y_{n-1}x_{n-1}^2+y_{n-1}^3-x_{n-1}y_{n-1}
\end{pmatrix},\\[13pt]
~~ x_0,y_0\in \mathbb{R},~ n\in \mathbb{N}_1.
\end{array}
\end{equation*}
\noindent The values of the parameters are chosen such that Lyapunov exponents are positive \cite{aa1} and therefore, the chaotic attractor with a spectacular dihedral symmetry, is presented in Fig. \ref{fig1} (a).
\noindent Another dihedral map \eqref{d1}, for $m=6$, with a spectacular chaotic attractor is presented in Appendix.

\subsection{Cyclic logistic map of IO}\label{sub22}

The map \eqref{d1} with $c\neq 0$, becomes the \textbf{cyclic logistic map}, with $C_m$ symmetry.

\begin{equation}\label{d111}
f(z)=(a+bz\bar{z}+ c i )z+d \bar{z}^{\;m-1}.
\end{equation}

\begin{theorem} The map \eqref{d111} presents cyclic symmetry.
\begin{proof}

Because, except the term $g(z)=iz$, all the other terms verify dihedral symmetry (see Theorem \ref{prop1}), one considers next only $g(z)$.
\begin{enumerate}
[leftmargin=0pt,align=left,itemindent=*,noitemsep]

\item
\emph{Reflection symmetry}

Due to the term $g(z)$, the reflection symmetry is not verified:
\[
S_0(g(z))=\overline{iz}=-i\bar{z}\neq i\bar{z}=g(S_0(z)),
\]
and, therefore, the reflection symmetry is broken.
\item
\vspace{5px}
\emph{Rotation symmetry}

With $iz=\mathrm{e}^{i\frac{\pi}{2}}r \mathrm{e}^{i\theta}$, one has
\[
R_1(g(z))=R_1(r\mathrm{e}^{i(\frac{\pi}{2}+\theta)})=\mathrm{e}^{i(\frac{\pi}{2}+\theta+\frac{2\pi}{m})}=\mathrm{e}^{i\frac{\pi}{2}}r\mathrm{e}^{i(\theta+\frac{2\pi}{m})}=iR_1(z)=g(R_1(z)).
\]
Therefore, the map \eqref{d111} maintains the cyclic symmetry, but looses the reflection symmetry.

The theorem is proved.
 \qedhere
\end{enumerate}
\end{proof}

\end{theorem}
\noindent Let $m=4$, i.e. the group $D_4$, for $a=-1.86$, $b=2.1$, $c=0.1$, $d=-1$. Then, the IVP of \eqref{d111}, in the cartesian computational form is

\begin{equation*}\label{f}
\begin{array}
[c]{ll}%
\begin{pmatrix}
x_n\\
y_n
\end{pmatrix}
=
\begin{pmatrix}
1.1x_{n-1}^3+5.1y_{n-1}^2x_{n-1}-1.86x_{n-1}-0.1y_{n-1}\\[4pt]
5.1x_{n-1}^2y_{n-1}+1.1y_{n-1}^3-1.86y_{n-1}+0.1x_{n-1}
\end{pmatrix},\\[13pt]~~x_0,y_0\in \mathbb{R},~n\in \mathbb{N}_1.
\end{array}
\end{equation*}
The chaotic attractor is presented in Fig. \ref{fig2} (a).
\noindent In order to exemplify computationally the symmetry  of this attractor, consider a point belonging to the attractor, e.g., $A(-0.085522,-0.9266)$, and calculate a reflection of this point with respect the $\Re(z)$ axis.

\noindent In this case the eight elements of $D_4$ are determined by the following rotation and reflection
\begin{equation*}
R_1=
\begin{pmatrix}
0& -1\\
1& 0
\end{pmatrix},
\end{equation*}
and
\[
S_0=
\begin{pmatrix}
1& 0\\
0& -1
\end{pmatrix},
\]
and the obtained point $B_1$ is given by
\[
B_1=S_0 A=
\begin{pmatrix}
 -0.0855\\[2pt]
 0.9266
 \end{pmatrix},
 \]
symmetric point with respect the horizontal axis but which, as can be seen from the image, does not belong to the attractor at an expected position.

\noindent On the other side, if one rotates the point $A$ with $R_1$, one obtains the point $B_2$
 \[
 B_2=R_1A=
 \begin{pmatrix}
 0.9266\\[2pt]
 -0.0855
 \end{pmatrix},
 \]
 i.e. a point which, as can be seen, belongs to the attractor.

\section{Symmetry breaking of maps \eqref{d1} of FO}
In this section it is proved that the dihedral and cyclic logistic maps \eqref{d1} lose their symmetry once they are transformed into FO maps via Caputo's discrete fractional difference.

\noindent Following \cite{xyz}, a convergent numerical method applied to an IVP, can be considered as a dynamical system whose dynamics reflect the behavior of the underlying system

Due to the fact that the FO variants of the IO dihedral and cyclic logistic maps loose symmetries, hereafter they will simply called logistic maps.

\subsection{Logistic map \eqref{d11} of FO}
 Let $N_a=\{a,a+1,a+2,\ldots\}$. Then, for~$q>0$ and $q\not\in \mathbb{N}$, the~$q$-th Caputo-like discrete fractional difference with starting point $a$ of a function $u:N_a\rightarrow \mathbb{R}$ is defined as~\cite{ati,xxx2,patru}
\begin{equation*}\label{capa}
\Delta_a^q u(t)=\Delta_a^{-(m-q)}\Delta^m u(t)=\frac{1}{\Gamma(m-q)}\sum_{s=a}^{t-(m-q)}(t-s-1)^{(m-q-1)}\Delta^m u(s),
\end{equation*}
for $t\in N_{a+m-q}$ and $m=[q]+1$.

\noindent $\Delta^m$ is the $m$-th order forward difference operator,
\[
\Delta^m u(s)=\sum_{k=0}^{m}\binom {n}{k}(-1)^{m-k}u(s+k),
\]
while $\Delta_a^{-q}$ represents the $q$-th fractional sum of $u$ starting from $a$,
\begin{equation*}\label{suma}
\Delta_a^{-q}u(t)=\frac{1}{\Gamma(q)}\sum_{s=a}^{t-q}(t-s-1)^{(q-1)}u(s),~t\in \mathbb{N}_{a+q},
\end{equation*}
and $t^{(q)}$ is the falling factorial, in the following form:
\begin{equation*}\label{fact}
t^{(q)}=\frac{\Gamma(t+1)}{\Gamma(t-q+1)}=t(t-1)...(t-q+1),
\end{equation*}
with the convention that at pole the division yields zero.

\noindent For the case $q\in(0,1)$, considered in this paper, $m=1$. Then, $\Delta^m u(s)=\Delta u(s)=u(s+1)-u(s)$, and for the usual case of starting point $a=0$, Caputo's fractional difference, denoted hereafter $\Delta_*^q$, becomes
\begin{equation*}
\Delta_*^q u(t)=\frac{1}{\Gamma(1-q)}\sum_{s=0}^{t-(1-q)}(t-s-1)^{(-q)}\Delta u(s).
\end{equation*}
Consider the following autonomous Caputo discrete IVP
\begin{equation}\label{eq2}
{\Delta_*^q u(t)=f(u(t+q-1)),~ t\in \mathbb{N}_{1-q}, ~u(0)=u_0,
}\end{equation}

\noindent {for $f:\mathbb{R}\rightarrow \mathbb{R}$ a scalar function and $u_0$ a real number.}

\begin{figure}
\begin{center}
\includegraphics[scale=0.75]{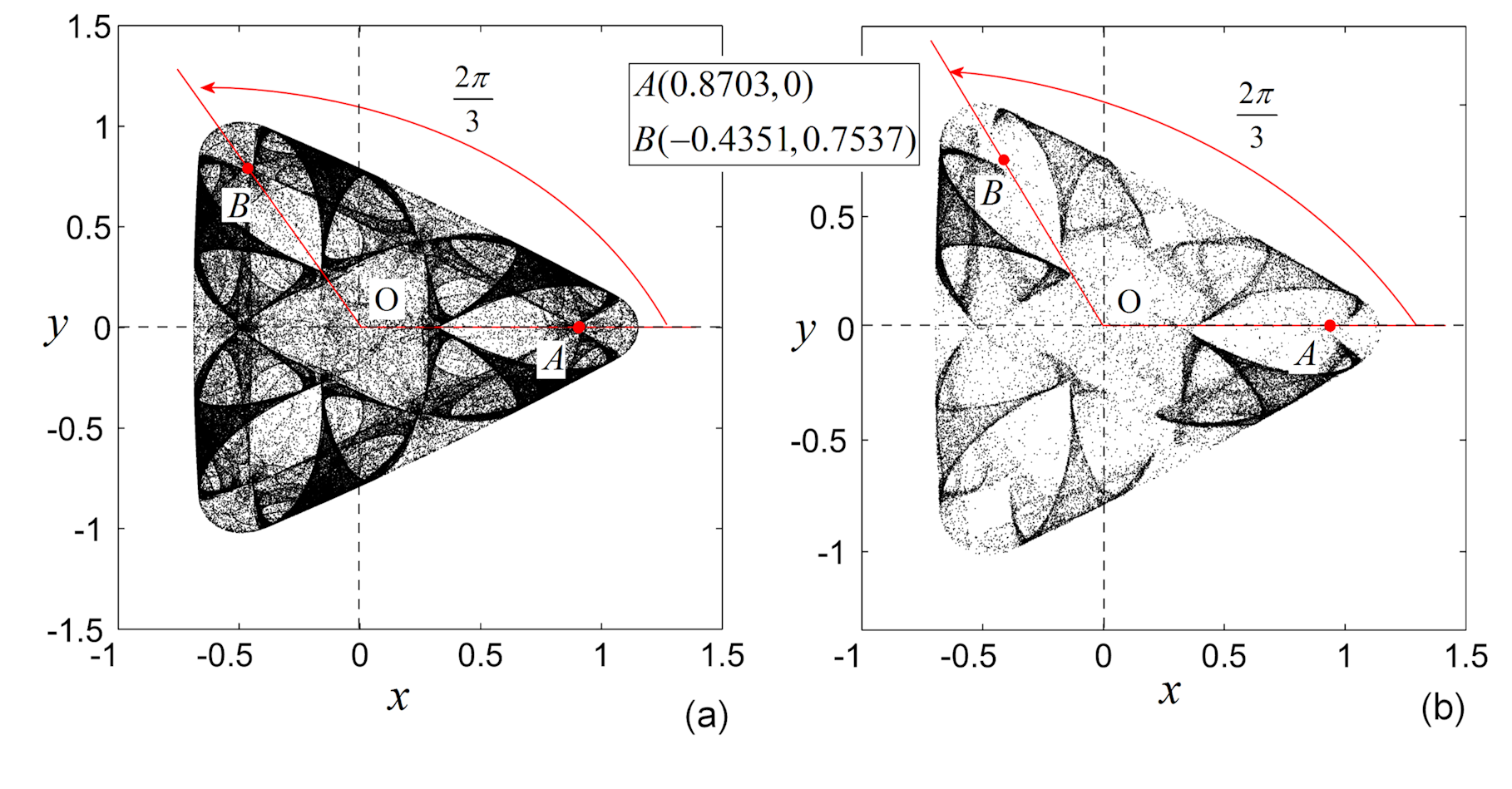}
\caption{Chaotic attractor of the logistic map \eqref{d11}. (a) IO case. (b) FO case. Point $B$ is obtained by rotating the point $A$ with $\frac{2\pi}{3}$.}
\label{fig1}
\end{center}
\end{figure}

\begin{figure}
\begin{center}
\includegraphics[scale=0.75]{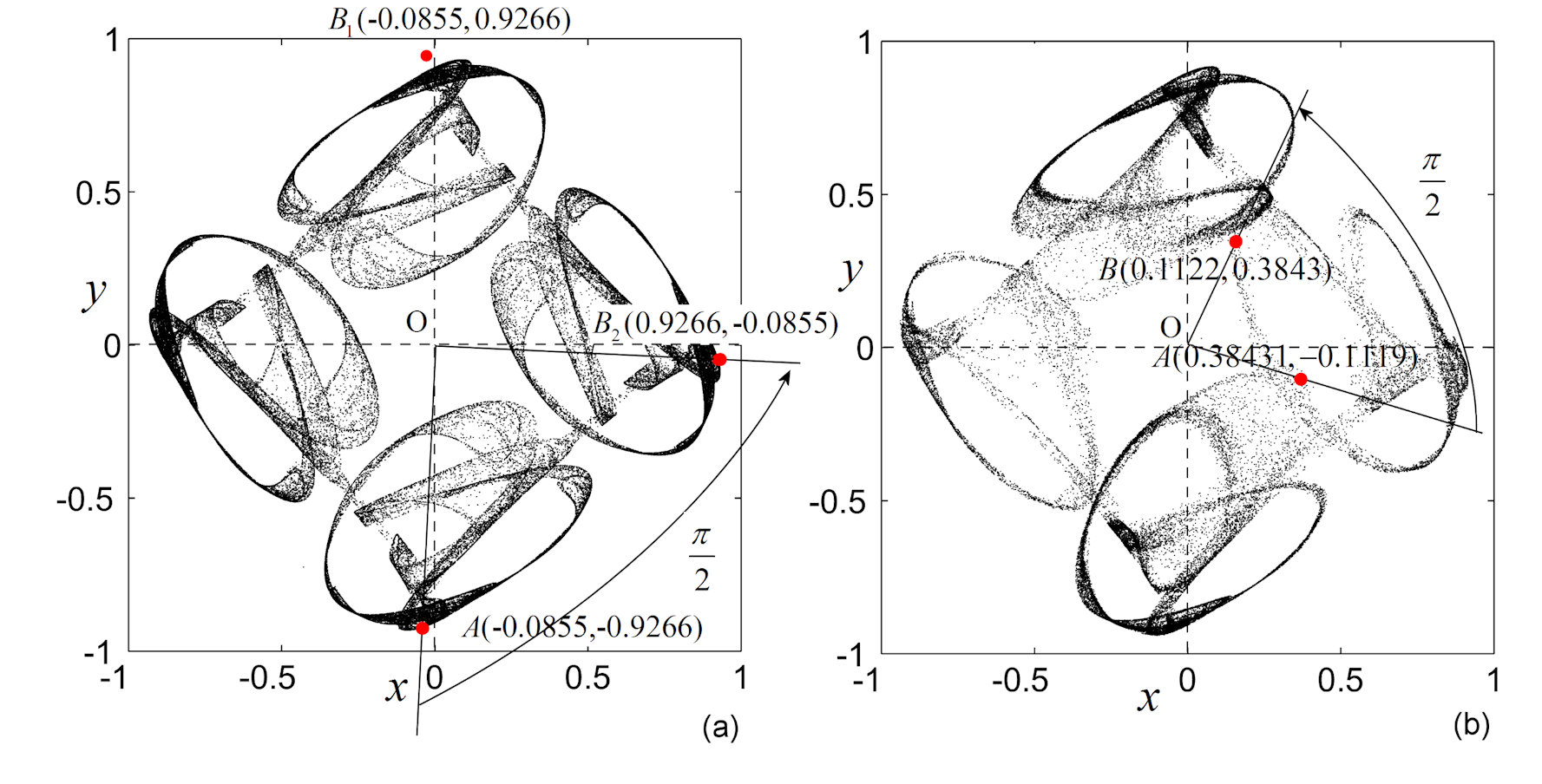}
\caption{Chaotic attractor of the logistic map \eqref{d111}. (a) IO case. (b) FO case. Left: Point $B_1$ is the reflection, with respect the horizontal axis, of the point $A$, while $B_2$ is image of the rotated point $A$ with $\frac{\pi}{2}$; Right: Point $B$ is the image of the rotated point $A$ with $\frac{\pi}{2}$.}
\label{fig2}
\end{center}
\end{figure}

\begin{theorem}\cite{xxx2}
$u:\mathbb{N}_1\rightarrow R$ is a solution of the IVP \eqref{eq2} if and only if $u(t)$ is solution to the following fractional Taylor difference formula
\begin{equation}\label{inte}
u(t)=u_0+\frac{1}{\Gamma(q)}\sum_{s=1-q}^{t-q}(t-s-1)^{(q-1)}f(u(s+q-1)).
\end{equation}
\end{theorem}
\noindent Existence and uniqueness results on the IVP \eqref{eq2} can be found in, e.g., \cite{ati,xxx2,trei}.
\noindent The convenable numerical form of \eqref{inte} is
\begin{equation}\label{eq3}
u(n)=u(0)+\frac{1}{\Gamma(q)}\sum_{k=1}^n\frac{\Gamma(n-k+q)}{\Gamma(n-k+1)}f(u(k-1)), ~{n\in \mathbb{N}_1~}.
\end{equation}
Consider next the IVP of FO \eqref{eq2} in the complex plane. With $u=z$, $z\in \mathbb{C}$ and $z(0)=x(0)+iy(0)$, the IVP has the following form

\begin{equation*}\label{ivi}
\Delta_*^q z(t)=f(z(t+q-1)),~~
 t\in \mathbb{N}_{1-q},~ z(0)=z_0.
\end{equation*}
Conform to \eqref{eq3}, the numerical integral is
\begin{equation}\label{eqq4}
z(n)=z(0)+\dfrac{1}{\Gamma(q)}\displaystyle\sum_{k=1}^n\dfrac{\Gamma(n-k+q)}{\Gamma(n-k+1)} f(z(k-1)),~~
n\in \mathbb{N}_1.
\end{equation}
With $f$ given by \eqref{d11}, the computational cartesian form becomes
\begin{equation}\label{eqqq}
\begin{array}
[c]{ll}
\begin{pmatrix}
x_n\\
y_n
\end{pmatrix}
=
\begin{pmatrix}
x_0+\dfrac{1}{\Gamma(q)}\displaystyle\sum_{k=1}^n\dfrac{\Gamma(n-k+q)}{\Gamma(n-k+1)}(-1.804x_{k-1}+x_{k-1}^3+x_{k-1}y_{k-1}^2+0.5~x_{k-1}^2-0.5~y_{k-1}^2)\\[4pt]
y_0+\dfrac{1}{\Gamma(q)}\displaystyle\sum_{k=1}^n\dfrac{\Gamma(n-k+q)}{\Gamma(n-k+1)}(-1.804~y_{k-1}+y_{k-1}x_{k-1}^2+y_{k-1}^3-x_{k-1}y_{k-1})
\end{pmatrix},
\\[6pt]
n\in \mathbb{N}_1.
\end{array}
\end{equation}
A chaotic attractor, for $q=0.03$ is presented in Fig. \ref{fig1} (b).

\begin{remark}\label{remus}
Note that the shape of the FO chaotic attractor approaches the shape IO chaotic attractor, for only small values $q\ll1$.
\end{remark}
\noindent The following result proves that for $y(0)\neq 0$, the fractionality in the Caputo's sense breaks the dihedral symmetry of the logistic map \eqref{d11}.

\begin{theorem}\label{th1}
For $q\in(0,1)$, and $y(0)\neq 0$, the logistic map \eqref{d11} of FO has no dihedral symmetry.
\begin{proof}

Denote for simplicity
\[
F_n(z)=\dfrac{1}{\Gamma(q)}\displaystyle\sum_{k=1}^n\dfrac{\Gamma(n-k+q)}{\Gamma(n-k+1)} f(z(k-1),
\]
so that \eqref{eqq4} becomes
\[
z(n)=z(0)+F_n(z), ~n\in \mathbb{N}_1.
\]
\begin{enumerate}
[leftmargin=0pt,align=left,itemindent=*,noitemsep]

\item
\emph{Reflection symmetry}

Because $S_0(z(n))=\bar{z}(n)$
\begin{equation}\label{x1}
f(S_0(z(n))=f(\bar{z}(n))=z(0)+F_n(\bar{z}(n)).
\end{equation}
On the other side
\begin{equation}\label{x2}
S_0(f(z(n))=\overline{z_0+F_n(z(n))}=\bar {z}(0)+\overline{F_n(z(n))}.
\end{equation}

Following the results in Subsection \ref{sub21}, it is easy to see that $\overline{F_n(z(n))}=F_n(\bar{z}(n))$. Moreover, if $y(0)\neq 0$, it follows that $z(0)=x(0)+iy(0)\neq x(0)-iy(0)=\bar{z}(0)$. Therefore, from relations \eqref{x1} and \eqref{x2}, one has $F_n(S_0(z(n)))\neq S_0(F_n(z(n)))$, and for $y(0)\neq0$, the logistic map \eqref{d11} of FO does not has reflection symmetry.
\item
\emph{Rotation symmetry}

\[
f(R_1(z(n)))=z(0)+F_n(R_1(z(n))).
\]
On the other side
\[
R_1(f(z(n)))=R_1(z(0))+R_1(F_n(z(n))).
\]
Like for the reflection property, previously (Subsection \ref{sub21}) it has been proved that $R_1(F_n(z(n)))=F_n(R_1(z(n)))$. However, because for $y(0)\neq0$, with $z(0)=\lvert z(0)\rvert\mathrm{e}^{i\frac{\pi}{2}}$, by comparing the last two relations, it follows that $R_1(z(0))=\lvert z(0)\rvert\mathrm{e}^{i(\frac{\pi}{2}+\frac{2\pi}{m})}\neq z(0)$.

\noindent Therefore, the logistic map of FO \eqref{d11} does not presents the rotation symmetry and, therefore, the map is not dihedral.

Theorem is proved.
 \qedhere
\end{enumerate}
\end{proof}
\end{theorem}
\noindent In Fig. \ref{fig1} (a) and (b) can be seen the difference between the rotation of some point, e.g., $A(0.8703,0)$, with $R_1$ in the cases of both IO and FO attractors. The obtained point, $B(-0.4351,0.7537)$, is given by
\[
B=R_1(A)=
\begin{pmatrix}
\cos \frac{2\pi}{3} -\sin\frac{2\pi}{3}\\[4pt]
 \sin\frac{2\pi}{3} \cos\frac{2\pi}{3}
 \end{pmatrix}
 \begin{pmatrix}
 0.8703\\[4pt]
 0
\end{pmatrix}
=
\begin{pmatrix}
-0.4351\\[4pt]
0.7537
\end{pmatrix}.
\]

\begin{figure}
\begin{center}
\includegraphics[scale=0.6]{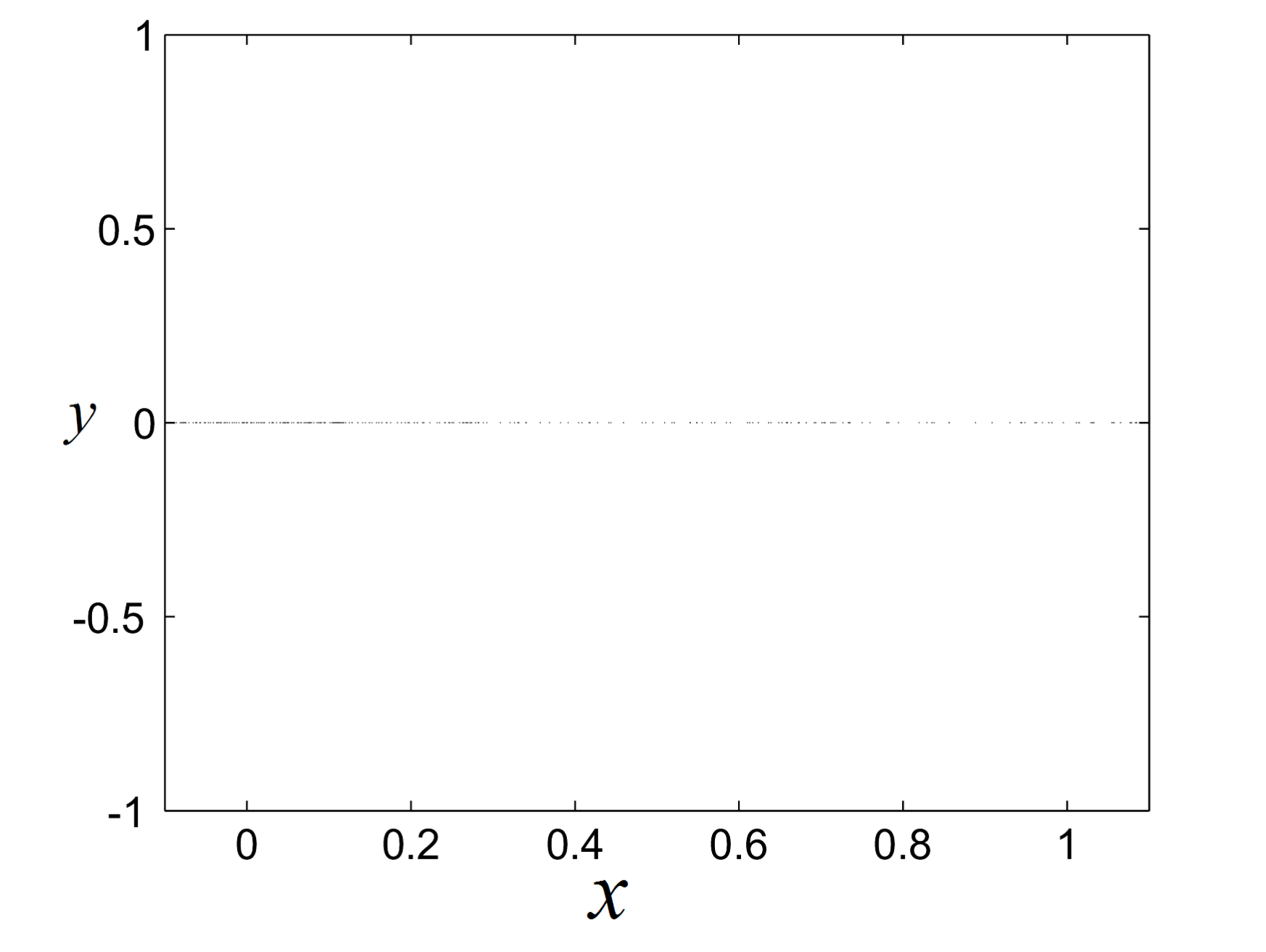}
\caption{One-dimensional dihedral chaotic attractor of the map \eqref{d11} of FO, for $y_0=0$.}
\label{fig3}
\end{center}
\end{figure}

\noindent As can be seen in the case of FO attractor (Fig. \ref{fig1} (b)), the rotation with $\frac{2\pi}{3}$ does not give the expected result.\\
\noindent For $y(0)=0$ one can see that the fractionality does not broke the dihedral symmetry. However, the case $y(0)=0$ and $x(0)\in \mathbb{R}$, generates a particular one-dimensional $D_1$ dihedral symmetry along the axis $\Re(z)$ (Fig. \ref{fig3}).\\
Theorem \ref{th1} can be extended to general continuous maps with or without symmetries.

\begin{theorem}
Let a map $f:\mathbb{C}\rightarrow \mathbb{C}$ continuous. For $q\in(0,1)$ and $y(0)\neq 0$, Caputo's FO variant of $f$, has no dihedral symmetry.
\end{theorem}
\noindent The continuity property is required to deduct the solution and to its existence.

\subsection{Logistic map \eqref{d111} of FO}
Let the FO variant of the cyclic logistic map \eqref{d111} and the underlying IVP
\[
\Delta_*^qz(t)=f(z(t+q-1)),~t\in \mathbb{N}_{1-q},~z(0)=z_0,
\]
where $f$ is given by \eqref{d111}.

\begin{theorem}\label{th2}
For $q\in (0,1)$ and $y(0)\neq 0$, the logistic map \eqref{d111} of FO has no cyclic symmetry.
\begin{proof}
See Theorem \ref{th1}.
\end{proof}
\end{theorem}
\noindent For $m=4$ and coefficients given in Section \ref{sub22}, the cartesian computational form of the solution is
\begin{equation*}
\begin{array}
[c]{ll}
\begin{pmatrix}
x_n\\
y_n
\end{pmatrix}
=
\begin{pmatrix}
x_0+\dfrac{1}{\Gamma(q)}\displaystyle\sum_{k=1}^n\dfrac{\Gamma(n-k+q)}{\Gamma(n-k+1)}(1.1x_{k-1}^3+5.1y_{k-1}^2x_{k-1}-1.86x_{k-1}-0.1y_{k-1})\\[4pt]
y_0+\dfrac{1}{\Gamma(q)}\displaystyle\sum_{k=1}^n\dfrac{\Gamma(n-k+q)}{\Gamma(n-k+1)}(5.1x_{k-1}^2y_{k-1}+1.1y_{k-1}^3-1.86y_{k-1}+0.1x_{k-1}
\end{pmatrix},\\[6pt]
~~n\in \mathbb{N}_1.
\end{array}
\end{equation*}
The chaotic non-symmetric attractor is shown in Fig. \ref{fig2} (b). The non-symmetry is exemplified by considering the point belonging to the attractor, $A(0.38431,-0.1119)$ which, after a rotation with $R_1$ of $\pi/2$, is moved at the point $B(0.1122,0.3843)$, which is not the symmetric of $A$.

\section{Bifurcation diagrams}
In this section, Bifurcation Diagrams (BDs), are considered for dihedral and cyclic logistic maps \eqref{d11} and \eqref{d111}, respectively, for both IO and FO case. For shortness, only the real component $x$ is considered.
In the IO case the diagrams are realized versus the parameter $a$, while in the FO case, diagrams are obtained versus the fractional order $q$.
\noindent The intersections between orthogonal lines through the BDs can be considered as Poincar\'{e}-like sections which, together with other tools, like time series or phase portrait, help to classify the attractors.

\subsection{IO variants}
\begin{figure}
\begin{center}
\includegraphics[scale=0.65]{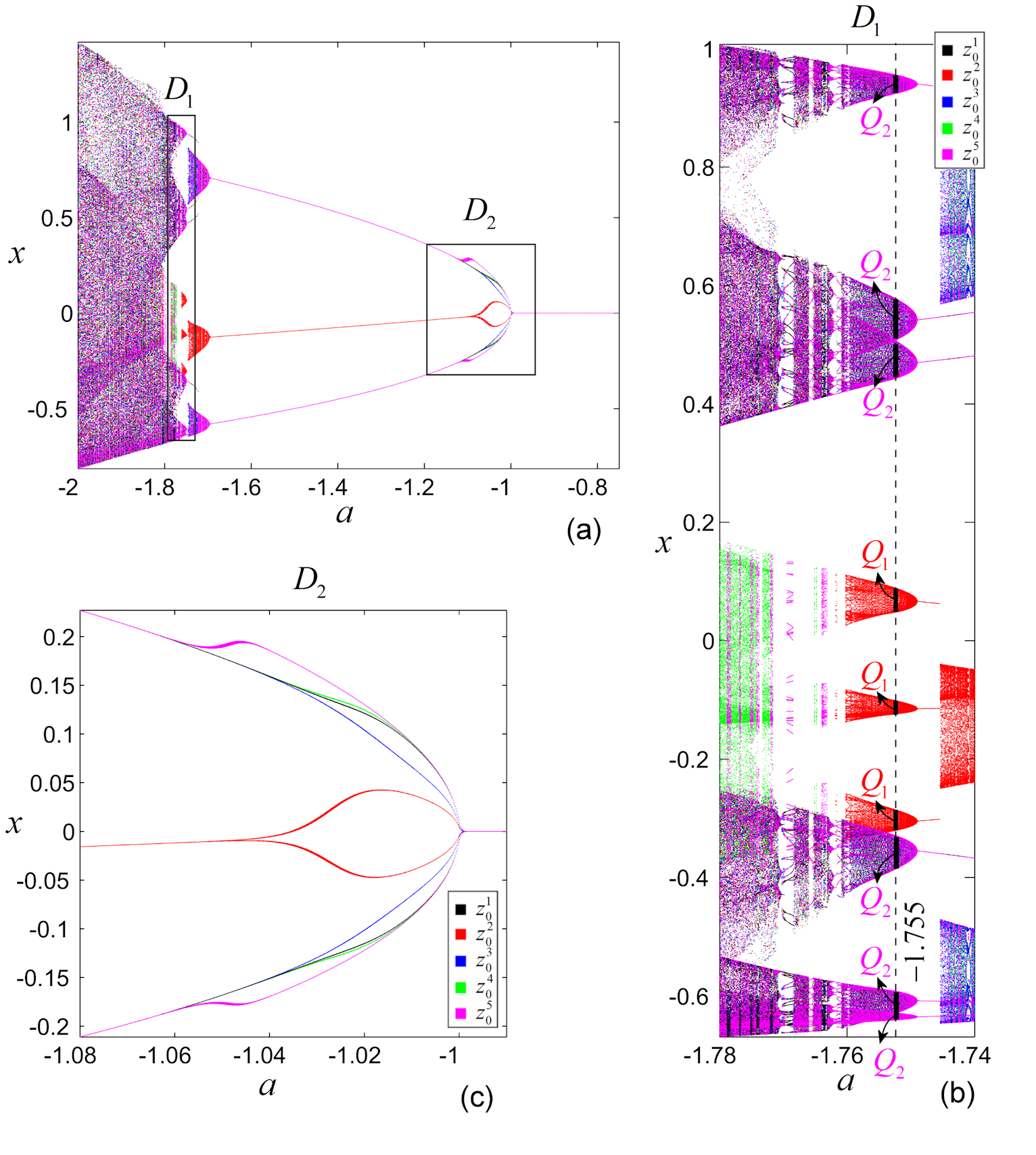}
\caption{Dihedral logistic map \eqref{d11} of IO. (a) Bifurcation diagram vs parameter $a$. (b) Zoomed detail $D_1$. Dotted line (Poincar\'{e}-like section) through $a=-1.755$ identifies the quasiperiodic attractors $Q_1$ and $Q_2$. (c) Zoomed detail $D_2$.}
\label{fig4}
\end{center}
\end{figure}

\begin{figure}
\begin{center}
\includegraphics[scale=0.75]{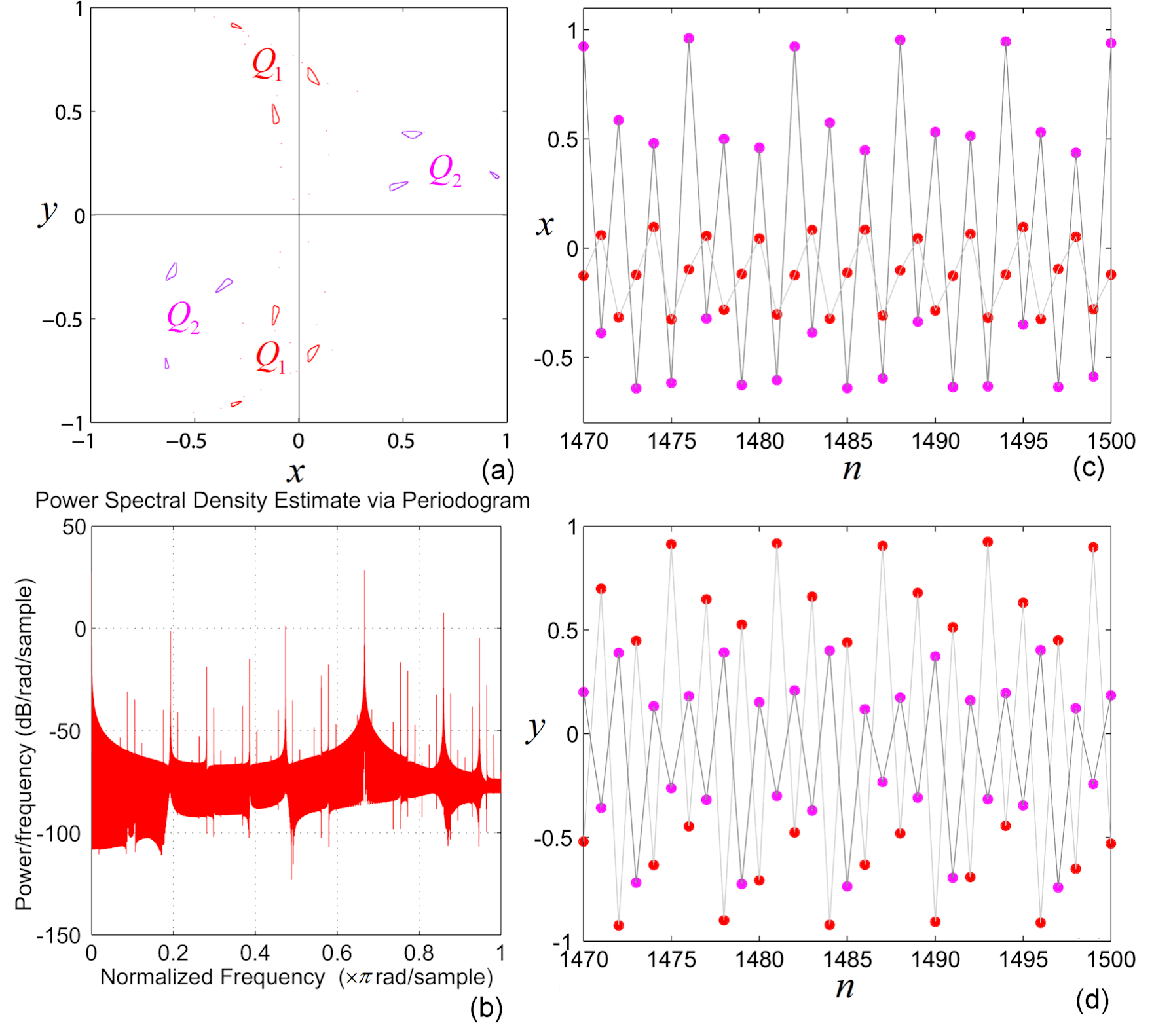}
\caption{Quasiperiodic attractors $Q_{1,2}$ (see Fig. \ref{fig4} (b)). (a) Phase plot. (b) PSD of the $x$ variable of $Q_1$. (c)-(d) Time series of both variables $x$ and $y$, respectively, revealing the visiting attractors order.}
\label{fig5}
\end{center}
\end{figure}

\begin{figure}
\begin{center}
\includegraphics[scale=0.75]{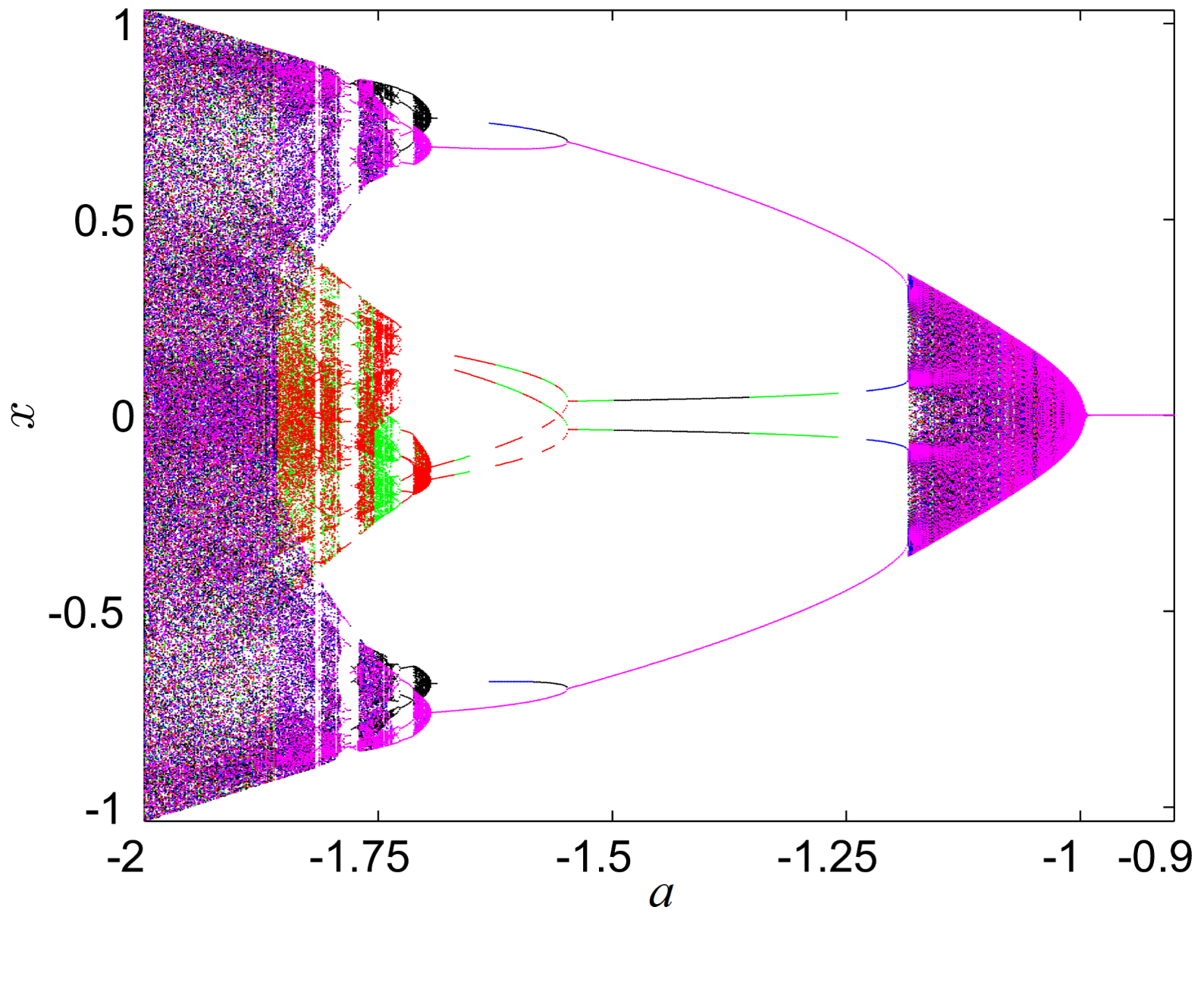}
\caption{Bifurcation diagram of the cyclic logistic map \eqref{d111} vs parameter $a$.}
\label{fig6}
\end{center}
\end{figure}

\begin{enumerate}
[leftmargin=0pt,align=left,itemindent=*,noitemsep]
\item
\emph{Dihedral logistic map }\eqref{d11}

The BD vs $a$ is presented in (Fig. \ref{fig4} (a)). Details $D_1$ and $D_2$, presented in Fig. \ref{fig4} (b) and (c), reveal the complexity of the system dynamics.

At $a=-1$, a bifurcation appears (detail $D_2$), whose analytical study is not considered in this paper.
\noindent

The diagram has been determined with five different empiric initial conditions $z_0^1=(x_0^1,y_0^1)=(0.05,0.1)$, $z_0^2=(x_0^2,y_0^2)=(0.01,0.01)$, $z_3^3=(x_0^3,y_0^3)=(0.001,0.9667)$, $z_0^4=(x_0^4,y_0^4)=-(0.477,-0.4965)$ and $z_0^5=(x_0^5,y_0^5)=(0.5,0.0001)$.
\noindent A useful notion in the analyze of the diagram is the Bifurcative Set (BS) introduced in \cite{plis}, i.e. the BD corresponding to a single initial condition.
\noindent For the considered initial conditions, there correspond BSs colored in red, blue, green, black and magenta, respectively. The BD could be considered as composed by all BSs. Even five initial conditions have been considered, the BD looks to be composed actually by only two BSs, red and magenta, while the other BSs identify with one or other of these two BSs.

\noindent Consider the case $a=-1.755$. As the Poincar\'{e}-like section through the BD (dotted line) shows, there exist two six-fold quasiperiodic attractors $Q_1$ and $Q_2$ (red and magenta, respectively) which are underlined by the tick black lines (Fig. \ref{fig4} (b)). The phase plot of these attractors is presented in Fig. \ref{fig5} (a), and the quasiperiodicity characteristic is revealed by the Power Spectrum Density (PSD) of the $x$ variable of $Q_1$ in Fig. \ref{fig5} (b), where beside the basic frequency and its modulation frequencies (the higher peaks), one can see the new born harmonics (the smaller peaks).
\noindent All the six branches of $Q_1$ and $Q_2$ can be remarked in the phase plot in Fig. \ref{fig5} (a) while the points visiting order, are highlighted by Fig. \ref{fig5} (c) and (d). To note that because of the symmetry, each attractors $Q_1$ and $Q_2$ are composed by two symmetric groups with three quasiperiodic islands.
\item
\vspace{5px}
\emph{Cyclic logistic map }\eqref{d111}

The bifurcation diagram vs $a$ is presented in Fig. \ref{fig6}. Again five initial conditions have been considered and, like for the dihedral map \eqref{d11}, only three different BSs appear. More initial conditions gave the same result, three BDs.\\

\end{enumerate}
\subsection{FO variants}
\begin{figure}
\begin{center}
\includegraphics[scale=0.9]{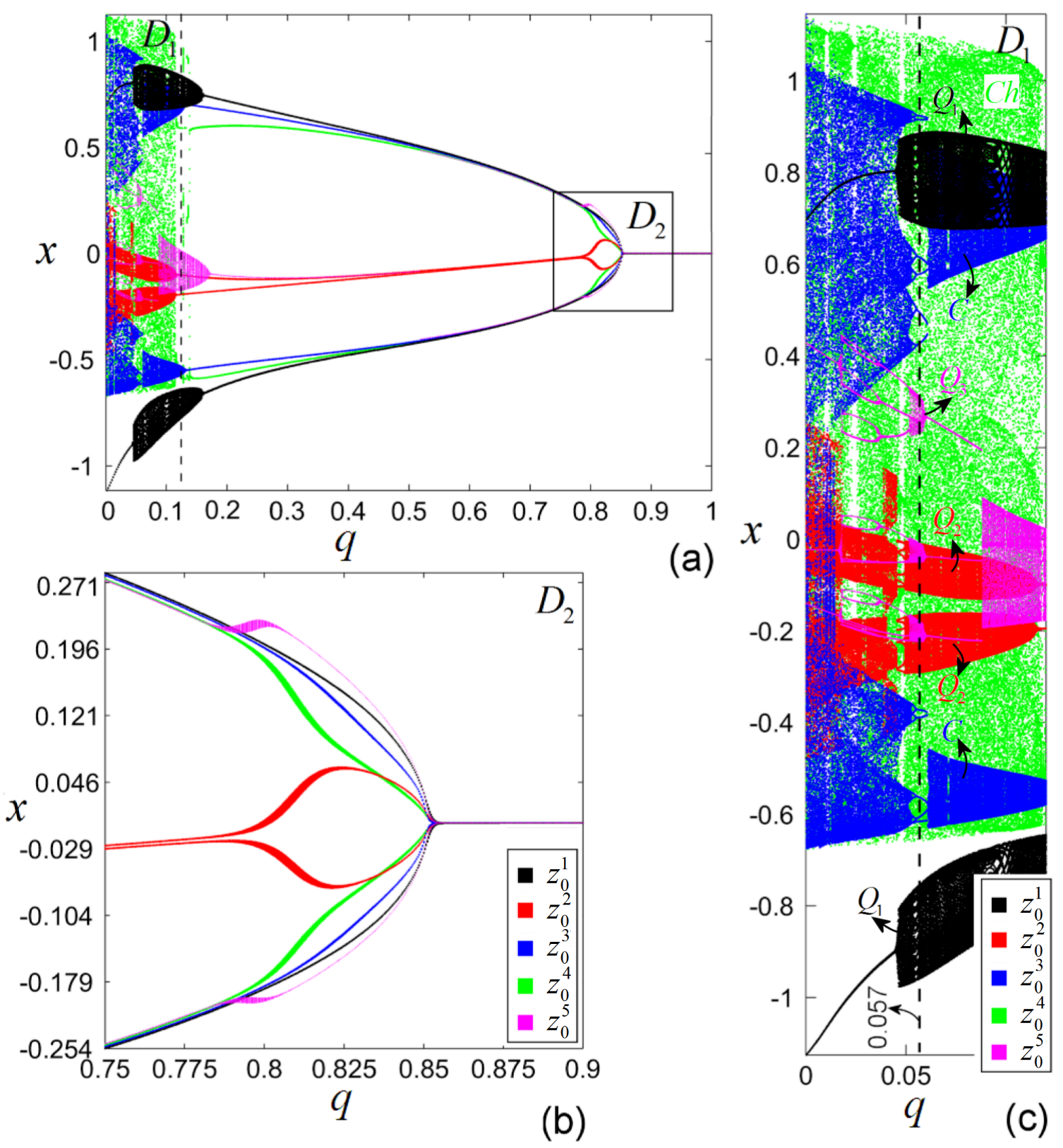}
\caption{Logistic map \eqref{d11} of FO. (a) Bifurcation diagram vs $q$, with five initial conditions. (b) Zoomed detail $D_2$. (c) Zoomed detail $D_1$. Dotted line (Poincar\'{e}-like section) through $q=0.057$, identifies five attractors: quasiperiodic-like ($Q_1$, $Q_2$ and $Q_3$), chaotic ($Ch$) and periodic-like attractor $C$. $z_0^i$, $i=1,2,...,5$ denote the initial conditions.}
\label{fig7}
\end{center}
\end{figure}

\begin{figure}
\begin{center}
\includegraphics[scale=0.85]{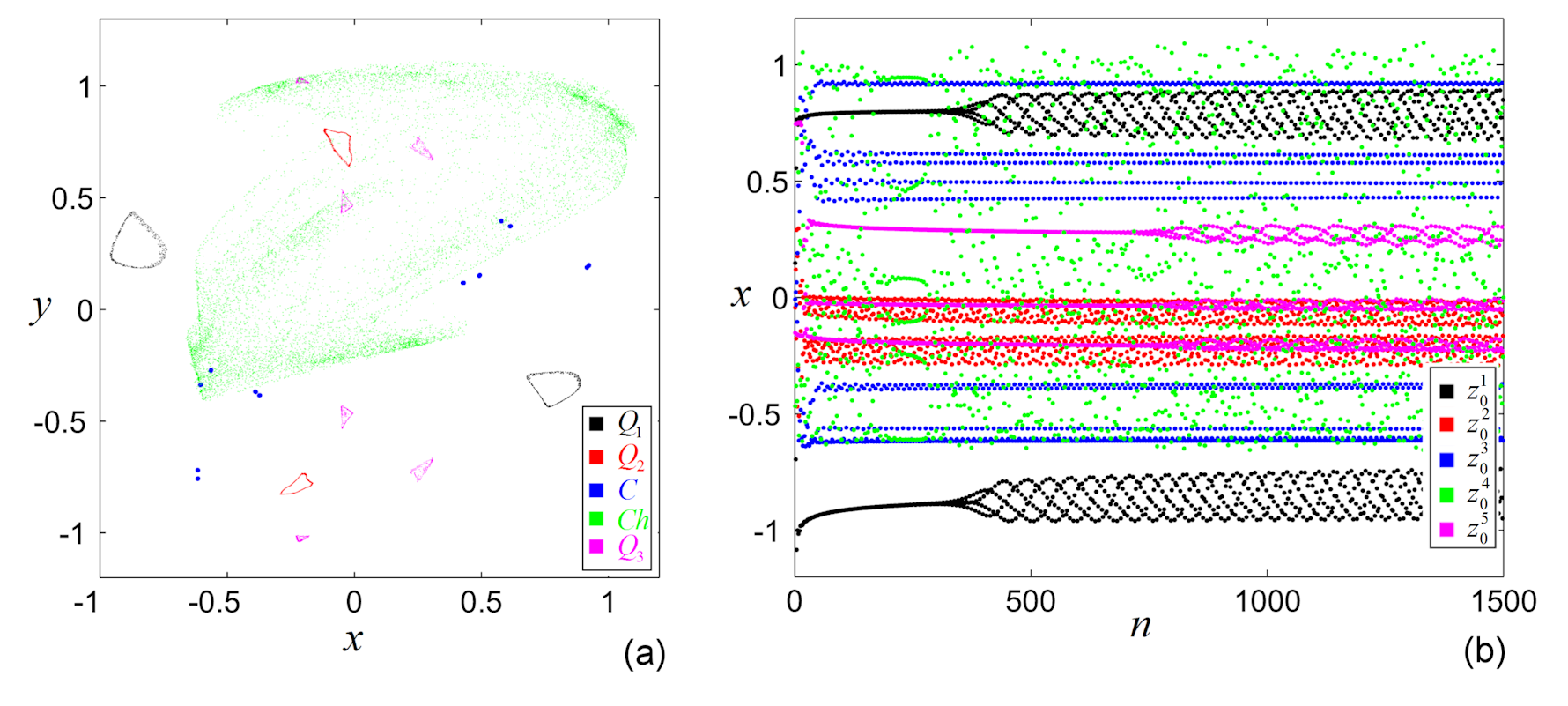}
\caption{Attractors of the map \eqref{d11} of FO for $q=0.057$ (see Fig. \ref{fig7} (c)). (a) Phase plot. (b) Time series.}
\label{fig8}
\end{center}
\end{figure}

\begin{figure}
\begin{center}
\includegraphics[scale=0.85]{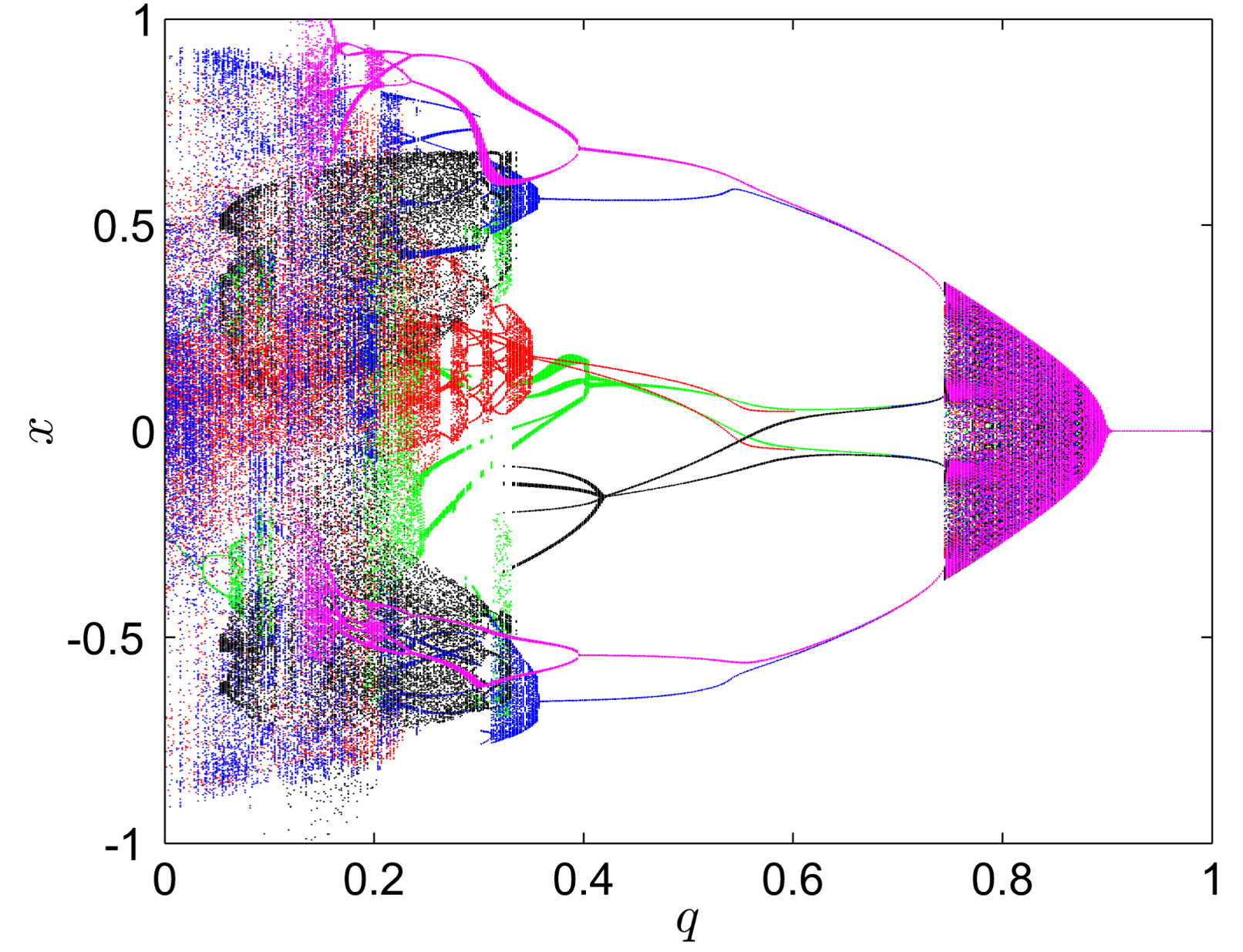}
\caption{Bifurcation diagram of the map \eqref{d111} vs $q$.}
\label{fig9}
\end{center}
\end{figure}

\begin{enumerate}
[leftmargin=0pt,align=left,itemindent=*,noitemsep]

\item

\emph{The logistic map \eqref{d11}}

The BD vs $q$ is presented in Fig. \ref{fig7} (a). Now, as can be seen, for $q=0.85$, each of the considered five initial conditions $z_0^1=(0.001,.9667)$, $z_0^2=(-0.477,-0.4965)$, $z_0^3=(0.5,.0001)$, $z_0^4=(-0.1,-0.1)$, $z_0^5=(0.00001,0.1)$, generates a different BS (black, red, blue, green and magenta, respectively). As for the IO counterpart, the FO variant preserves the interesting behavior of the BD of IO within a neighborhood of $a=-1$ (see the zoom around $q=0.85$ in Fig. \ref{fig7} (b)). However, while in the case of the dihedral logistic map of IO, it is natural to consider the bifurcation phenomenon, the situation changes dramatically in the case of continuous-time or discrete-time systems of FO where, except some particular cases, such as S-asymptotical systems of FO \cite{nepe} or delayed FO systems \cite{nepe2}, periodicity in FO systems, is not allowed \cite{nepe,neper1}. Beside the analytical proofs in the mentioned references, an intuitive reasoning of this phenomenon is the time history of FO orbits or trajectories. Therefore, the bifurcation notion in the case of FO systems has no meaning.
In this paper, the orbit which apparently looks periodic or quasiperiodic will be called periodic-like or quasiperiodic-like, respectively.\\
The detail $D_2$ presents five different attractors: three quasiperiodic-like attractors $Q_{1,2,3}$ (black, red, and magenta, respectively), a chaotic attractor $Ch$ (green), and a periodic-like orbit (cycle-like), $C$ (blue). The attractors corresponding to the values $q=0.057$ (dotted line) are depicted in Fig. \ref{fig8} (a) and (b)  (phase plot and time series, respectively). Because the considered value of $q$ ($0.057$) is relatively distanced from $0$ (see Remark \ref{remus}), where the FO chaotic attractor resembles with its symmetric counterpart of the IO case, in this case the chaotic attractor (green) looks strongly deformed.
\vspace{5px}
\item

\emph{The logistic map \eqref{d111}}

The BD vs $q$ is presented in Fig. \ref{fig9}. As in the case of the map \eqref{d11}, the five initial conditions generate five different BSs.

\end{enumerate}

\begin{figure}
\begin{center}
\includegraphics[scale=1]{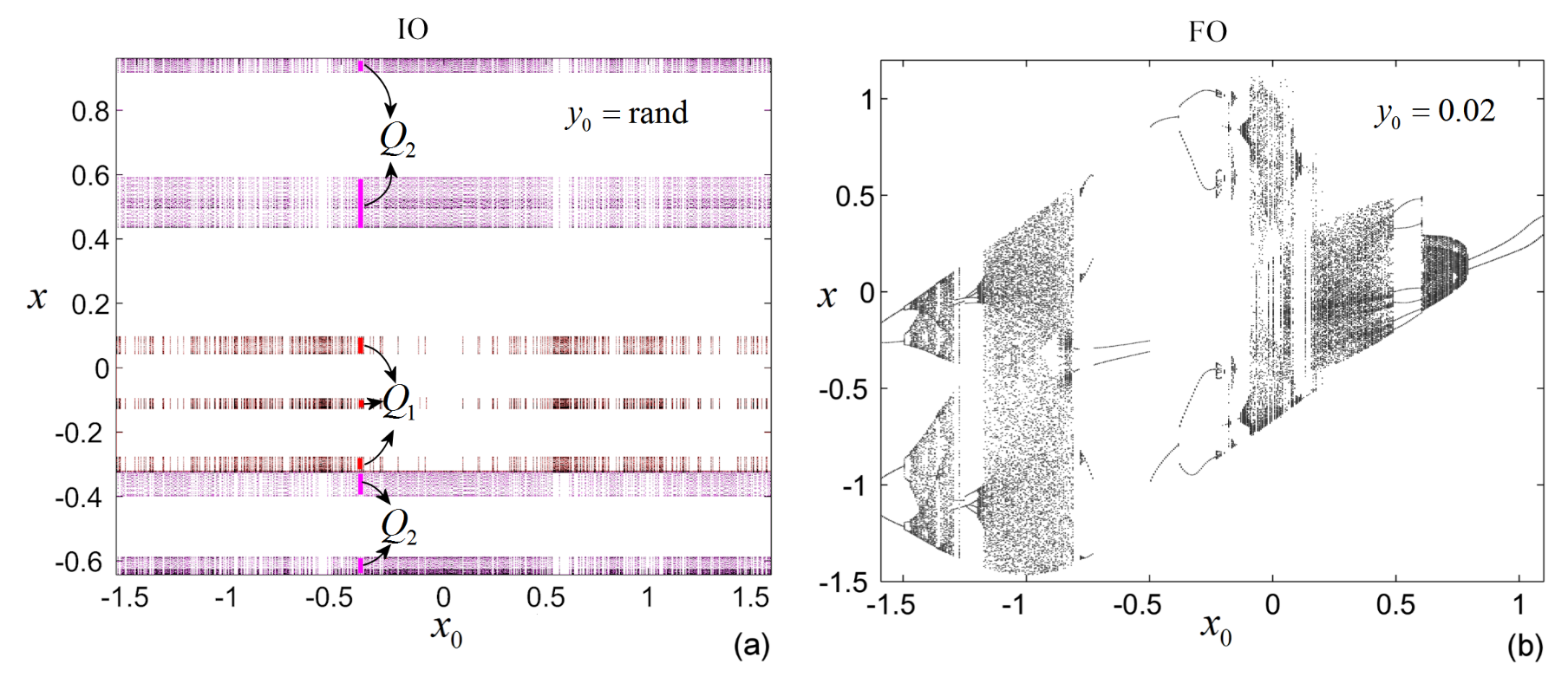}
\caption{Bifurcation diagrams of the map \eqref{d11} vs the first component of initial condition, $x_0$. (a) IO case where, beside $x0$ variable within the range $[-1.5,1.5]$, $y_0$ is changed randomly every bifurcation step. For every $x_0$ and random $y_0$, one obtain the same two quasiperiodic attractors $Q_{1,2}$. (b) FO case where $y_0$ is fixed.}
\label{fig10}
\end{center}
\end{figure}

\noindent To underline the difference between the IO and the FO cases, consider the dihedral logistic map \eqref{d11}, and the logistic map \eqref{eqqq} for $a=1.755$ and $q=0.057$, respectively, in the following numerical experiment: a BD vs one of the two components of the initial condition, namely $x_0$, is generated for both cases. In the IO case, every step in the BD vs $x_0$, the second component $y_0$ is also chosen randomly in order to verify the result in a stronger variation of the initial conditions. In the FO case, $y_0$ is maintained constant. As Fig. \ref{fig10} shows, to each $x_0$ and $y_0$ variable, one obtains the same two quasiperiodic attractors $Q_{1,2}$, obtained also in Fig. \ref{fig4} (b). In the FO case, as the BD shows, for each value of $x_0$ one obtains another attractor, fact which underlines the conclusion that, the BD vs the fractional order is strongly dependent on initial conditions. The result is similar for the cyclic map.

\begin{remark}
\begin{enumerate}[noitemsep,labelindent=0pt, leftmargin=*]\quad

\item [ i)]
Since the statement that the BD in the case of discrete FO systems it is not yet rigorously proved, this dependence on initial conditions can be conjuctured as holding true in other discrete and even continuous-time FO systems with Caputo's approach \cite{v2,v3vduma,plis};
\item[ii)] The dependence on initial conditions seems to be significant for values of $q$ less than about $0.5$ (see also \cite{v2},\cite{v3},\cite{duma});
\item[iii)] Contrary to intuition, if for a set of parameters a IO map presents  a chaotic evolution, the BD of its FO counterpart shows that once that $q\rightarrow1$, the FO system behaves regularly-like, not chaotic as expected.

\end{enumerate}
\end{remark}

\section*{Conclusion}
In this paper it is proved that the discrete FO variant of discrete maps in the complex plane, obtained via Caputo's discrete fractional differences, generally cannot admit dihedral symmetry. The example of logistic dihedral map is considered. Despite the symmetries of the considered IO maps, except a particular case of the initial conditions, the FO variants cannot have symmetry. The results extend to general maps, showing that a map with or without symmetry, looses its symmetry once it is transformed into a FO variant, in the Caputo sense. Another important property in the case of FO discrete systems is computationally revealed: the notion of bifurcation diagram versus the FO cannot be rigourously determined as for the IO case where one or two initial conditions are sufficient. It is conjuctured that in FO systems, the bifurcation diagrams versus the FO depends strongly on the initial conditions: every initial condition leads to a different diagram.
This phenomenon, which also appeared in other discrete and even continuous-time systems, seems to be typical for about the first half of the range $(0,1)$ of the fractional order.\\

\noindent \textbf{Funding}: No funding to declare.

\noindent \textbf{Declaration of Competing Interest}: Author declare that he does has no conflict of interest.

\noindent \textbf{Acknowledgment~} Author thanks Alexandru-David Abrudan, master student, Babes-Bolyai University, Cluj-Napoca.

 \section*{APPENDIX}\label{app}

\begin{figure}
\begin{center}
\includegraphics[scale=1]{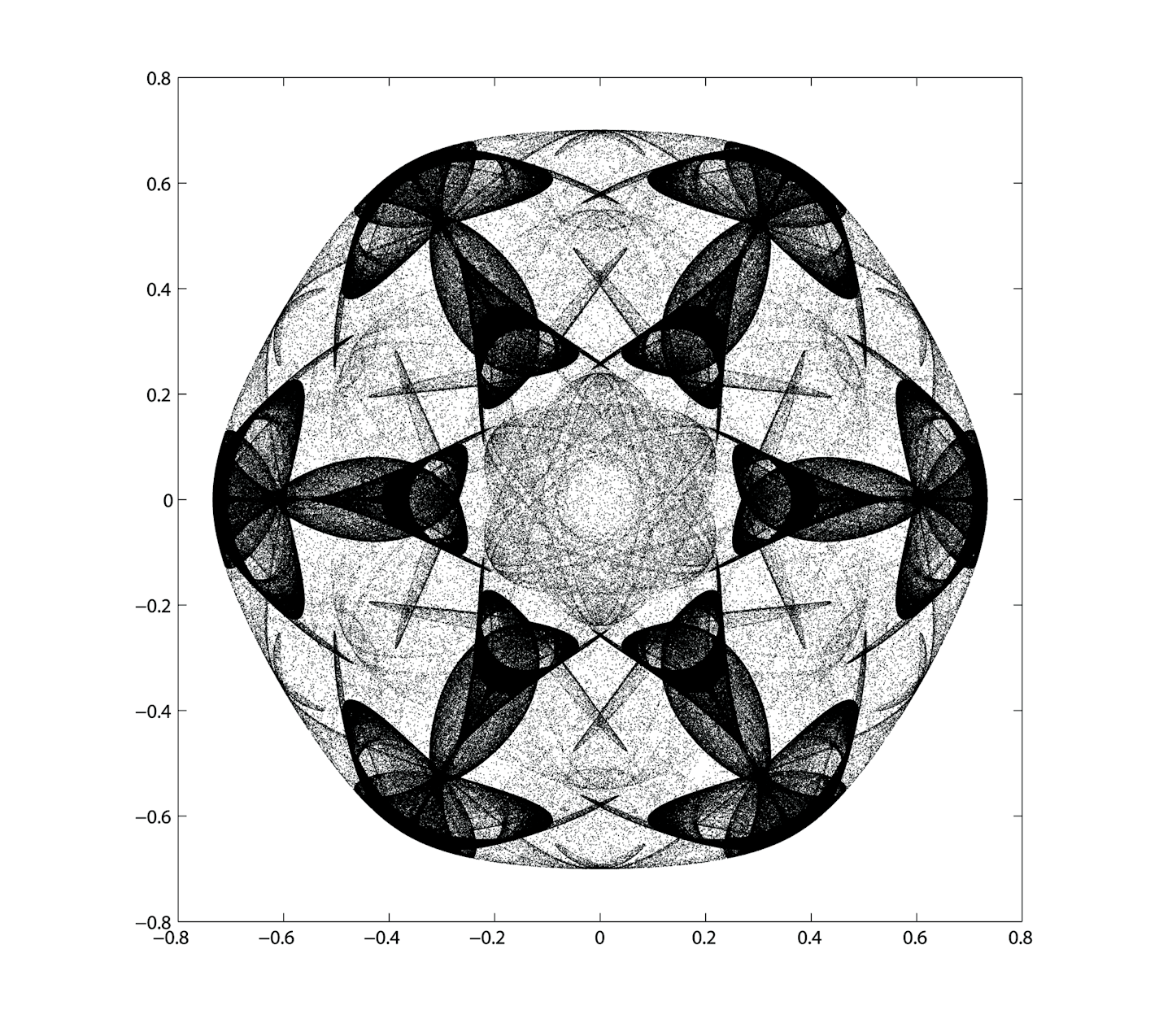}
\caption{Dihedral logistic map \eqref{sase} for $m=6$.  }
\label{fig11}
\end{center}
\end{figure}

Another dihedral logistic map \eqref{d1}, for $m=6$, is

\begin{equation}\label{sase}
f(z)=(\alpha+\beta z\bar{z}+\gamma \Re(z^n))z+\delta \bar{z}^{\; m-1},
\end{equation}
with $\alpha=-2.584$, $\beta=5$, $\gamma=-2$ and $\delta=-1$. $\Re(z^n)$ can be computed as follows: $\Re(z^n)=(z^n+\bar{z}^n)/2$.
The real and imaginary components are

\begin{equation*}
\begin{array}
[c]{ll}
f_1(x,y)=-2.5840x+5x^3+5xy^2-2x^7+30x^5y^2-30x^3y^4+2y^6x-x^5+10x^3y^2-5xy^4,\\[4pt]
f_2(x,y)=2y^7+30x^4y^3-2x^6y+5x^2y+5y^3-10x^2y^3+y^5+5x^4y-2.5840y-30x^2y^5,
\end{array}
\end{equation*}
\noindent and the chaotic attractor is presented in Fig. \ref{fig11}.

\newpage{\pagestyle{empty}\cleardoublepage}

\newpage{\pagestyle{empty}\cleardoublepage}

\newpage{\pagestyle{empty}\cleardoublepage}
\newpage{\pagestyle{empty}\cleardoublepage}

\end{document}